\documentstyle[amstex]{amsart}
\input psfig.sty

\newtheorem{lemma}{Lemma}
\newtheorem{cor}{Corollary}
\newtheorem{prop}{Proposition}
 
\newtheorem{theorem}{Theorem}

\theoremstyle{definition}
\newtheorem{defn}{Definition}

\theoremstyle{remark}

\newcommand{\Q}{{\Bbb Q}}

\newcommand{\Gal}{\operatorname{Gal}}
\newcommand{\LCM}{\operatorname{LCM}}

\begin{document}

\title[Diagonals of a Regular Polygon]{The Number
	of Intersection Points Made by the Diagonals of a Regular Polygon}
\date{November 18, 1997}

\author{Bjorn Poonen}
\thanks{The first author is supported by an NSF Mathematical
Sciences Postdoctoral Research Fellowship.
Part of this work was done at MSRI,
where research is supported in part by NSF grant DMS-9022140.}
\address{AT\&T Bell Laboratories \\ Murray Hill, NJ 07974, USA}
\curraddr{University of California at Berkeley \\ Berkeley, CA 94720-3840, USA}
\email{poonen@@math.berkeley.edu}

\author{Michael Rubinstein}
\address{AT\&T Bell Laboratories, Murray Hill, NJ 07974, USA}
\curraddr{University of Waterloo \\ Waterloo, ON N2L 3G1, Canada}
\email{mrubinst@@uwaterloo.ca}

\begin{abstract}
We give a formula for the number of interior intersection
points made by the diagonals of a regular $n$-gon.
The answer is a polynomial on each residue class modulo 2520.
We also compute the number of regions formed by the diagonals,
by using Euler's formula $V-E+F=2$.
\end{abstract}

\maketitle

%****************************************************************************
\section{Introduction}
\label{intro}

\begin{figure}[p]
\centerline{\psfig{file=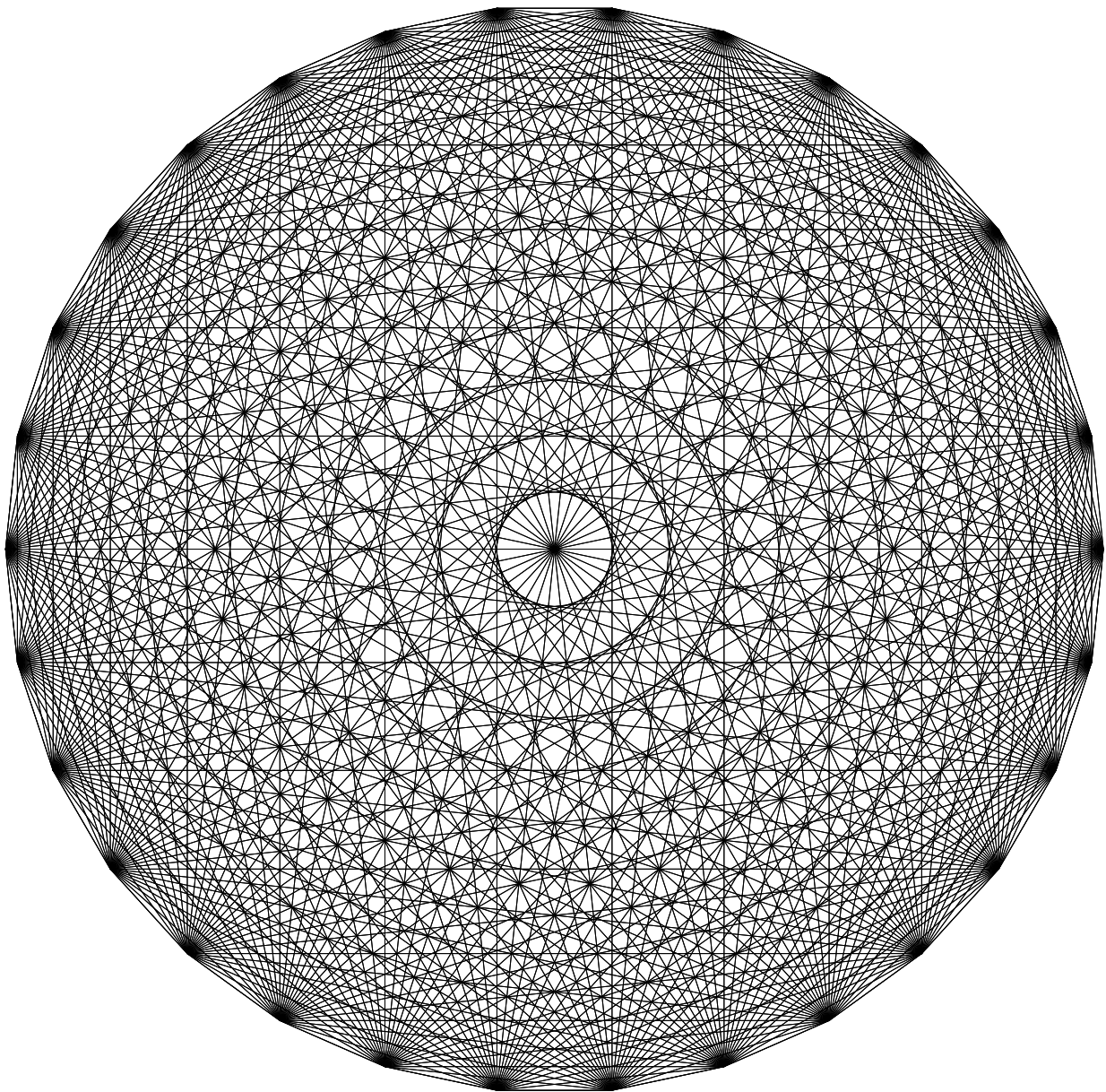,width=6.5in}}
\caption{ The 30-gon with its diagonals.
There are 16801 interior intersection
points: 13800 two line intersections, 2250 three line intersections,
420 four line intersections, 180 five line intersections, 120 six
line intersections, 30 seven line intersections, and 1 fifteen
line intersection.}
\label{30gon}
\end{figure}

We will find a formula for the number $I(n)$ of intersection points
formed inside a regular $n$-gon by its diagonals.
The case $n=30$ is depicted in Figure~\ref{30gon}.
For a {\em generic} convex $n$-gon, the answer would be $n \choose 4$,
because every four vertices would be the endpoints of a unique pair
of intersecting diagonals.
But $I(n)$ can be less, because in a regular $n$-gon
it may happen that three or more diagonals meet at an interior point,
and then some of the $n \choose 4$ intersection points will coincide.
In fact, if $n$ is even and at least~6,
$I(n)$ will always be less than $n \choose 4$,
because there will be $n/2 \ge 3$ diagonals meeting at the center point.
It will result from our analysis that for $n>4$,
the maximum number of diagonals of the regular $n$-gon
that meet at a point other than the center is
$$\begin{array}{rl}
	2	& \text{if $n$ is odd},	\\
	3	& \text{if $n$ is even but not divisible by 6},	\\
	5	& \text{if $n$ is divisible by 6 but not 30, and},	\\
	7	& \text{if $n$ is divisible by 30}.
\end{array}$$
with two exceptions: this number is~2 if $n=6$, and~4 if $n=12$.
In particular, it is impossible to have~8 or more diagonals
of a regular $n$-gon meeting at a point other than the center.
Also, by our earlier remarks, the fact that no three diagonals
meet when $n$ is odd will imply that $I(n)={n \choose 4}$ for odd $n$.

A careful analysis of the possible configurations of three diagonals
meeting will provide enough information to permit us in theory
to deduce a formula for $I(n)$.
But because the explicit description of these configurations is so complex,
our strategy will be instead to use this information
to deduce only the {\em form} of the answer,
and then to compute the answer for enough small $n$
that we can determine the result precisely.
The computations are done in Mathematica, Maple and C,
and annotated source codes can be obtained via anonymous ftp
at {\verb+http://math.berkeley.edu/~poonen+}.

In order to write the answer in a reasonable form, we define
	$$\delta_m(n) = \begin{cases}
				1	& \text{if $n \equiv 0 \pmod m$,} \\
				0	& \text{otherwise.}
			\end{cases}$$

\begin{theorem}
\label{countintersections}
For $n \ge 3$,
\begin{eqnarray*}
	I(n)	& = &	{n \choose 4}
			+ (-5 n^3 + 45 n^2 - 70 n + 24)/24 \cdot \delta_2(n)
			- (3n/2) \cdot \delta_4(n) \\
		&   &	\mbox{} + (-45 n^2 + 262n)/6 \cdot \delta_6(n)
			+ 42n \cdot \delta_{12}(n)
			+ 60 n \cdot \delta_{18}(n)	\\
		&   &	\mbox{} + 35n \cdot \delta_{24}(n)
			- 38n \cdot \delta_{30}(n)
			- 82n \cdot \delta_{42}(n)
			- 330n \cdot \delta_{60}(n)	\\
		&   &	\mbox{} - 144n \cdot \delta_{84}(n)
			- 96n \cdot \delta_{90}(n)
			- 144n \cdot \delta_{120}(n)
			- 96n \cdot \delta_{210}(n).
\end{eqnarray*}
\end{theorem}

Further analysis, involving Euler's formula $V-E+F=2$,
will yield a formula for the number $R(n)$ of regions
that the diagonals cut the $n$-gon into.

\begin{theorem}
\label{countregions}
For $n \ge 3$,
\begin{eqnarray*}
	R(n)	& = &	(n^4 - 6 n^3 + 23 n^2 - 42 n + 24)/24	\\
		&   &	\mbox{} 
			+ (-5 n^3 + 42 n^2 - 40 n - 48)/48 \cdot \delta_2(n)
			- (3n/4) \cdot \delta_4(n) \\
		&   &	\mbox{} + (-53 n^2 + 310n)/12 \cdot \delta_6(n)
			+ (49n/2) \cdot \delta_{12}(n)
			+ 32 n \cdot \delta_{18}(n)	\\
		&   &	\mbox{} + 19n \cdot \delta_{24}(n)
			- 36n \cdot \delta_{30}(n)
			- 50n \cdot \delta_{42}(n)
			- 190n \cdot \delta_{60}(n)	\\
		&   &	\mbox{} - 78n \cdot \delta_{84}(n)
			- 48n \cdot \delta_{90}(n)
			- 78n \cdot \delta_{120}(n)
			- 48n \cdot \delta_{210}(n).
\end{eqnarray*}
\end{theorem}

These problems have been studied by many authors before,
but this is apparently the first time the
correct formulas have been obtained.
The Dutch mathematician Gerrit Bol~\cite{bol}
gave a complete solution in~1936,
except that a few of the coefficients in his formulas are wrong.
(A few misprints and omissions in Bol's paper
are mentioned in~\cite{rigby2}.)

The approaches used by us and Bol are similar in many ways.
One difference (which is not too substantial)
is that we work as much as possible with roots of unity
whereas Bol tended to use more trigonometry
(integer relations between sines of rational multiples of $\pi$).
Also, we relegate much of the work to the computer,
whereas Bol had to enumerate the many cases by hand.
The task is so formidable that it is amazing to us
that Bol was able to complete it,
and at the same time not so surprising that it would contain a few errors!

Bol's work was largely forgotten.
In fact, even we were not aware of his paper
until after deriving the formulas ourselves.
Many other authors in the interim solved special cases of the problem.
Steinhaus~\cite{steinhausproblem} posed the problem
of showing that no three diagonals meet internally when $n$ is prime,
and this was solved by Croft and Fowler~\cite{croft}.
(Steinhaus also mentions this in~\cite{steinhaus},
which includes a picture of the 23-gon and its diagonals.)
In the 1960s, Heineken~\cite{heineken} gave a delightful argument
which generalized this to all odd $n$,
and later he~\cite{heineken2} and Harborth~\cite{harborth}
independently enumerated all three-diagonal intersections
for $n$ not divisible by 6.

The classification of three-diagonal intersections
also solves Colin Tripp's problem~\cite{tripp}
of enumerating ``adventitious quadrilaterals,''
those convex quadrilaterals for which the angles
formed by sides and diagonals are all rational multiples of $\pi$.
See Rigby's paper~\cite{rigby2} or the summary~\cite{rigby} for details.
Rigby, who was aware of Bol's work, mentions that Monsky and Pleasants
also each independently classified all three-diagonal intersections
of regular $n$-gons.
Rigby's papers partially solve Tripp's further problem of proving
the existence of all adventitious quadrangles using only
elementary geometry; i.e., without resorting to trigonometry.

All the questions so far have been in the Euclidean plane.
What happens if we count the interior intersections made by the diagonals
of a hyperbolic regular $n$-gon?
The answers are exactly the same, as pointed out in~\cite{rigby2},
because if we use Beltrami's representation of points of the
hyperbolic plane by points inside a circle in the Euclidean plane,
we can assume that the center of the hyperbolic $n$-gon corresponds
to the center of the circle, and then the hyperbolic $n$-gon with its
diagonals looks in the model exactly like a Euclidean regular $n$-gon
with its diagonals.
It is equally easy to see that the answers will be the same in
elliptic geometry.

%****************************************************************************
\section{When do three diagonals meet?}
\label{threemeet}

\begin{figure}
\centerline{\psfig{file=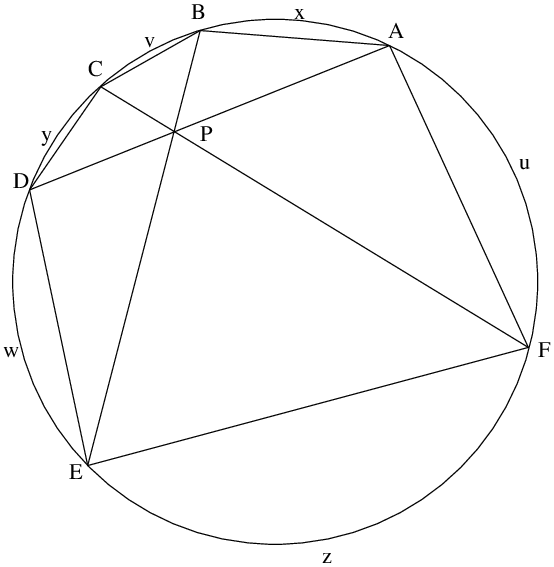,width=3in}}
\caption{\vrule depth 20pt width0pt}
\label{geometry}
\end{figure}

We now begin our derivations of the formulas for $I(n)$ and $R(n)$.
The first step will be to find a criterion for the concurrency
of three diagonals.
Let $A,B,C,D,E,F$ be six distinct points in order on
a unit circle dividing up the circumference into arc lengths
$u,x,v,y,w,z$  and assume that the three chords $AD, BE, CF$
meet at $P$ (see Figure~\ref{geometry}).

By similar triangles, $AF/CD = PF/PD$, 
$BC/EF = PB/PF$, 
$DE/AB = PD/PB$. 
Multiplying these together yields 
\begin{displaymath}
(AF \cdot BC \cdot DE)
/(CD \cdot EF \cdot AB)
= 1,
\end{displaymath}
and so
\begin{equation}
\label{product}
\sin(u/2) \sin(v/2) \sin(w/2) 
= \sin(x/2) \sin(y/2) \sin(z/2).
\end{equation}

Conversely, suppose six distinct points $A,B,C,D,E,F$
partition the circumference of a unit circle into 
arc lengths $u,x,v,y,w,z$ and suppose that~(\ref{product})
holds. Then the three diagonals 
$AD, BE, CF$ meet in a single point which we see as follows.
Let lines $AD$ and $BE$ intersect at $P_0$ .
Form the line through $F$ and $P_0$
and let $C'$ be the other intersection point of $FP_0$ 
with the circle. This partitions the circumference into
arc lengths $u,x,v',y',w,z$. As shown above, we have
\begin{displaymath}
\sin(u/2) \sin(v'/2) \sin(w/2)
= \sin(x/2) \sin(y'/2) \sin(z/2)
\end{displaymath}
and since we are assuming that~(\ref{product}) holds
for $u,x,v,y,w,z$ we get
\begin{displaymath}
	\frac{\sin(v'/2)}{\sin(y'/2)} = \frac{\sin(v/2)}{\sin(y/2)}.
\end{displaymath}
Let $\alpha=v + y = v' + y'$.  Substituting $v=\alpha-y$,
$v'=\alpha-y'$ above we get
\begin{displaymath}
	\frac{\sin(\alpha/2) \cos(y'/2)
	- \cos(\alpha/2) \sin(y'/2)}{\sin(y'/2)}
      = \frac{\sin(\alpha/2) \cos(y/2)
	- \cos(\alpha/2) \sin(y/2)}{\sin(y/2)}
\end{displaymath}
and so
\begin{displaymath}
	\cot(y'/2) = \cot(y/2).
\end{displaymath}
Now $0<\alpha/2<\pi$, so $y=y'$ and hence $C=C'$. Thus, 
the three diagonals $AD, BE, CF$ meet at a single point. 

So~(\ref{product}) gives a necessary and sufficient
condition (in terms of arc lengths) for the chords $AD, BE, CF$
formed by six distinct points $A,B,C,D,E,F$
on a unit circle to meet at a single point. 
In other words, to give an explicit answer to the question
in the section title, we need to characterize the positive
rational solutions to
\begin{eqnarray}
\label{product2}
	\sin(\pi U) \sin(\pi V) \sin(\pi W)
		& = & \sin(\pi X) \sin(\pi Y) \sin(\pi Z)	\\
\nonumber	U+V+W+X+Y+Z & = & 1.
\end{eqnarray}
(Here $U=u/(2\pi)$, etc.)
This is a trigonometric diophantine equation in the sense of~\cite{conway},
where it is shown that in theory, there is a finite computation
which reduces the solution of such equations to
ordinary diophantine equations.
The solutions to the analogous equation with only two sines
on each side are listed in~\cite{myerson}.

If in~(\ref{product2}), we substitute
$\sin(\theta)=(e^{i \theta}-e^{-i \theta})/(2i)$,
multiply both sides by $(2i)^3$, and expand,
we get a sum of eight terms on the left equalling
a similar sum on the right, but two terms on the left
cancel with two terms on the right since $U+V+W=1-(X+Y+Z)$, leaving
\begin{eqnarray*}
\lefteqn{-e^{i\pi (V+W-U)} + e^{-i\pi(V+W-U)}
	- e^{i\pi(W+U-V)} + e^{-i\pi(W+U-V)}
	- e^{i\pi(U+V-W)} + e^{-i\pi(U+V-W)} = }  \\
&& -e^{i\pi (Y+Z-X)} + e^{-i\pi(Y+Z-X)} - e^{i\pi(Z+X-Y)}
	+ e^{-i\pi(Z+X-Y)} - e^{i\pi(X+Y-Z)} + e^{-i\pi(X+Y-Z)}.
\end{eqnarray*}
If we move all terms to the left hand side,
convert minus signs into $e^{-i \pi}$, multiply by $i=e^{i\pi/2}$, and let
\begin{eqnarray*}
	\alpha_1 & = & V+W-U-1/2	\\
	\alpha_2 & = & W+U-V-1/2	\\
	\alpha_3 & = & U+V-W-1/2	\\
	\alpha_4 & = & Y+Z-X+1/2	\\
	\alpha_5 & = & Z+X-Y+1/2	\\
	\alpha_6 & = & X+Y-Z+1/2,
\end{eqnarray*}
we obtain
\begin{equation}
\label{thetwelve}
	\sum_{j=1}^6 e^{i \pi \alpha_j}
	+ \sum_{j=1}^6 e^{-i \pi \alpha_j} = 0,
\end{equation}
in which $\sum_{j=1}^6 \alpha_j = U+V+W+X+Y+Z = 1$.
Conversely, given rational numbers
$\alpha_1,\alpha_2,\alpha_3,\alpha_4,\alpha_5,\alpha_6$
(not necessarily positive) which sum to 1 and satisfy~(\ref{thetwelve}),
we can recover $U,V,W,X,Y,Z$, (for example, $U=(\alpha_2+\alpha_3)/2+1/2$),
but we must check that they turn out positive.

%****************************************************************************
\section{Zero as a sum of 12 roots of unity}
\label{twelve}

In order to enumerate the solutions to~(\ref{product2}), we are led,
as in the end of the last section, to classify the ways
in which 12 roots of unity can sum to zero.
More generally, we will study relations of the form
\begin{equation}
\label{relation}
	\sum_{i=1}^k a_i \eta_i = 0,
\end{equation}
where the $a_i$ are positive integers, and the $\eta_i$
are distinct roots of unity.
(These have been studied previously by Schoenberg~\cite{schoenberg},
Mann~\cite{mann}, Conway and Jones~\cite{conway}, and others.)
We call $w(S)=\sum_{i=1}^k a_i$ the {\em weight} of the relation $S$.
(So we shall be particularly interested in relations of weight~12.)
We shall say the relation~(\ref{relation}) is {\em minimal}
if it has no nontrivial subrelation; i.e., if
	$$\sum_{i=1}^k b_i \eta_i = 0, \quad a_i \ge b_i \ge 0$$
implies either $b_i=a_i$ for all $i$ or $b_i=0$ for all $i$.
By induction on the weight, any relation can be represented
as a sum of minimal relations (but the representation need not be unique).

Let us give some examples of minimal relations.
For each $n \ge 1$, let $\zeta_n=\exp(2 \pi i/n)$
be the standard primitive $n$-th root of unity.
For each prime $p$, let $R_p$ be the relation
	$$1 + \zeta_p + \zeta_p^2 + \cdots + \zeta_p^{p-1} = 0.$$
Its minimality follows from the irreducibility of the cyclotomic polynomial.
Also we can ``rotate'' any relation by multiplying
through by an arbitrary root of unity to obtain a new relation.
In fact, Schoenberg~\cite{schoenberg} proved that every relation
(even those with possibly negative coefficients)
can be obtained as a linear combination with positive
and negative integral coefficients of the $R_p$ and their rotations.
But we are only allowing positive combinations,
so it is not clear that these are enough to generate all relations.

In fact it is not even true!
In other words, there are other minimal relations.
If we subtract $R_3$ from $R_5$, cancel the 1's and
incorporate the minus signs into the roots of unity,
we obtain a new relation
\begin{equation}
\label{r5r3}
	\zeta_6 + \zeta_6^{-1} + \zeta_5
	+ \zeta_5^2 + \zeta_5^3 + \zeta_5^4 = 0,
\end{equation}
which we will denote $(R_5:R_3)$.
In general, if $S$ and $T_1,T_2,\ldots,T_j$ are relations,
we will use the notation $(S:T_1,T_2,\ldots,T_j)$
to denote any relation obtained by rotating the $T_i$
so that each shares exactly one root of unity with $S$
which is different for each $i$, subtracting them from $S$,
and incorporating the minus signs into the roots of unity.
For notational convenience, we will write $(R_5:4R_3)$
for $(R_5:R_3,R_3,R_3,R_3)$, for example.
Note that although $(R_5:R_3)$ denotes unambiguously (up to rotation)
the relation listed in~(\ref{r5r3}), in general there will be
many relations of type $(S:T_1,T_2,\ldots,T_j)$ up to rotational equivalence.
Let us also remark that including $R_2$'s in the list of $T$'s has no effect.

It turns out that recursive use of the construction above
is enough to generate all minimal relations of weight up to 12.
These are listed in Table~\ref{minimals}.
The completeness and correctness of the table will be
proved in Theorem~\ref{table} below.
Although there are 107 minimal relations up to rotational equivalence,
often the minimal relations within one of our classes are Galois conjugates.
For example, the two minimal relations of type $(R_5:2R_3)$ are
conjugate under $\Gal(\Q(\zeta_{15})/\Q)$, as pointed out in~\cite{mann}.

The minimal relations with $k \le 7$ ($k$ defined as in~(\ref{relation}))
had been previously catalogued in~\cite{mann},
and those with $k \le 9$ in~\cite{conway}.
In fact, the $a_i$ in these never exceed~1,
so these also have weight less than or equal to~9.

\begin{table}
\begin{center}
\begin{tabular}{c|c|c|}
	Weight	& Relation type &
			Number of relations of that type \\ \hline \hline
	2	& $R_2$	& 1	\\ \hline
	3	& $R_3$	& 1	\\ \hline
	5	& $R_5$	& 1	\\ \hline
	6	& $(R_5:R_3)$	& 1	\\ \hline
	7	& $(R_5:2R_3)$	& 2	\\ \cline{2-3}
		& $R_7$	& 1	\\ \hline
	8	& $(R_5:3R_3)$	& 2	\\ \cline{2-3}
		& $(R_7:R_3)$	& 1	\\ \hline
	9	& $(R_5:4R_3)$	& 1	\\ \cline{2-3}
		& $(R_7:2R_3)$	& 3	\\ \hline
	10	& $(R_7:3R_3)$	& 5	\\ \cline{2-3}
		& $(R_7:R_5)$	& 1	\\ \hline
	11	& $(R_7:4R_3)$	& 5	\\ \cline{2-3}
		& $(R_7:R_5,R_3)$	& 6	\\ \cline{2-3}
		& $(R_7:(R_5:R_3))$	& 6	\\ \cline{2-3}
		& $R_{11}$	& 1	\\ \hline
	12	& $(R_7:5R_3)$	& 3	\\ \cline{2-3}
		& $(R_7:R_5,2R_3)$	& 15	\\ \cline{2-3}
		& $(R_7:(R_5:R_3),R_3)$	& 36	\\ \cline{2-3}
		& $(R_7:(R_5:2R_3))$	& 14	\\ \cline{2-3}
		& $(R_{11}:R_3)$	& 1	\\ \hline
\end{tabular}
\end{center}
\caption{The 107 minimal relations of weight up to 12.}
\label{minimals}
\end{table}

\begin{theorem}
\label{table}
Table~\ref{minimals} is a complete listing of the minimal
relations of weight up to 12 (up to rotation).
\end{theorem}

The following three lemmas will be needed in the proof.

\begin{lemma}
\label{mann1}
If the relation~(\ref{relation}) is minimal,
then there are distinct primes $p_1<p_2<\cdots<p_s \le k$
so that each $\eta_i$ is a $p_1 p_2 \cdots p_s$-th root of unity,
after the relation has been suitably rotated.
\end{lemma}

\begin{pf}
This is a corollary of Theorem~1 in~\cite{mann}.
\end{pf}

\begin{lemma}
\label{2p}
The only minimal relations (up to rotation) involving only
the $2p$-th roots of unity, for $p$ prime, are $R_2$ and $R_p$.
\end{lemma}

\begin{pf}
Any $2p$-th root of unity is of the form $\pm \zeta^i$.
If both $+\zeta^i$ and $-\zeta^i$ occurred in the same relation,
then $R_2$ occurs as a subrelation.
So the relation has the form
	$$\sum_{i=0}^{p-1} c_i \zeta_p^i = 0$$
By the irreducibility of the cyclotomic polynomial,
$\{1,\zeta_p,\ldots,\zeta_p^{p-1}\}$ are independent over $\Q$
save for the relation that their sum is zero, so all the $c_i$ must be equal.
If they are all positive, then $R_p$ occurs as a subrelation.
If they are all negative, then $R_p$ rotated by -1 (i.e., 180 degrees)
occurs as a subrelation.
\end{pf}

\begin{lemma}
\label{pigeonhole}
Suppose $S$ is a minimal relation, and $p_1<p_2<\cdots<p_s$
are picked as in Lemma~\ref{mann1} with $p_1=2$ and $p_s$ minimal.
If $w(S)<2 p_s$, then $S$ (or a rotation) is of the
form $(R_{p_s}:T_1,T_2,\ldots,T_j)$ where the $T_i$ are minimal relations
not equal to $R_2$ and involving only $p_1 p_2 \cdots p_{s-1}$-th
roots of unity, such that $j<p_s$ and
	$$\sum_{i=1}^j [w(T_i)-2] = w(S) - p_s.$$
\end{lemma}

\begin{pf}
Since every $p_1 p_2 \ldots p_s$-th root of unity is uniquely expressible
as the product of a $p_1 p_2 \ldots p_{s-1}$-th root of unity
and a $p_s$-th root of unity, the relation can be rewritten as
\begin{equation}
\label{peel}
	\sum_{i=0}^{p_s-1} f_i \zeta_{p_s}^i = 0,
\end{equation}
where each $f_i$ is a sum of $p_1 p_2 \ldots p_{s-1}$-th roots of unity,
which we will think of as a sum (not just its value).

Let $K_m$ be the field obtained by adjoining the $p_1 p_2 \ldots p_m$-th
roots of unity to $\Q$.
Since
$[K_s:K_{s-1}]=\phi(p_1 p_2 \cdots p_s)/\phi(p_1 p_2 \cdots p_{s-1})
=\phi(p_s)=p_s-1$,
the only linear relation satisfied
by $1,\zeta_{p_s},\ldots,\zeta_{p_s}^{p_s-1}$ over $K_{s-1}$
is that their sum is zero.
Hence~(\ref{peel}) forces the values of the $f_i$ to be equal.

The total number of roots of unity in all the $f_i$'s is $w(S)<2 p_s$,
so by the pigeonhole principle, some $f_i$ is zero or
consists of a single root of unity.
In the former case, each $f_j$ sums to zero,
but at least two of these sums contain at least one root of unity,
since otherwise $s$ was not minimal,
so one of these sums gives a subrelation of $S$, contradicting its minimality.
So some $f_i$ consists of a single root of unity.
By rotation, we may assume $f_0=1$.
Then each $f_i$ sums to 1, and if it is not simply the single root
of unity~1, the negatives of the roots of unity in $f_i$ together
with~1 form a relation $T_i$ which is not $R_2$ and involves
only $p_1 p_2 \cdots p_{s-1}$-th roots of unity, and it is clear
that $S$ is of type $(R_{p_s}:T_{i_1},T_{i_2},\ldots,T_{i_j})$.
If one of the $T$'s were not minimal, then it could be decomposed
into two nontrivial subrelations, one of which would not share a root
of unity with the $R_{p_s}$, and this would give a nontrivial subrelation
of $S$, contradicting the minimality of $S$.
Finally, $w(S)$ must equal the sum of the weights of $R_{p_s}$
and the $T$'s, minus $2j$ to account for the roots of unity that
are cancelled in the construction
of $(R_{p_s}:T_{i_1},T_{i_2},\ldots,T_{i_j})$.
\end{pf}

\begin{pf*}{Proof of Theorem~\ref{table}}
We will content ourselves with proving that every relation of
weight up to 12 can be decomposed into a sum of the ones listed
in Table~\ref{minimals}, it then being straightforward to check that
the entries in the table are distinct, and that none of them can be
further decomposed into relations higher up in the table.

Let $S$ be a minimal relation with $w(S) \le 12$.
Pick $p_1<p_2<\cdots<p_s$ as in Lemma~\ref{mann1} with $p_1=2$
and $p_s$ minimal.
In particular, $p_s \le 12$, so $p_s=2$,3,5,7, or $11$.

\medskip
\noindent{\em Case 1: } $p_s \le 3$

Here the only minimal relations are $R_2$ and $R_3$, by Lemma~\ref{2p}.

\medskip
\noindent{\em Case 2: } $p_s=5$

If $w(S)<10$, then we may apply Lemma~\ref{pigeonhole} to
deduce that $S$ is of type
	$(R_5:T_1,T_2,\ldots,T_j)$
Each $T$ must be $R_3$ (since $p_{s-1} \le 3$), and $j=w(S)-5$ by
the last equation in Lemma~\ref{pigeonhole}.
The number of relations of type $(R_5:jR_3)$, up to rotation,
is ${5 \choose j}/5$.
(There are $5 \choose j$ ways to place the $R_3$'s, but one must
divide by 5 to avoid counting rotations of the same relation.)

If $10 \le w(S) \le 12$, then write $S$ as in~(\ref{peel}).
If some $f_i$ consists of zero or one roots of unity, then the
argument of Lemma~\ref{pigeonhole} applies, and $S$ must be of the
form $(R_5:jR_3)$ with $j \le 4$,
which contradicts the last equation in the Lemma.
Otherwise the numbers of (sixth) roots of unity occurring
in $f_0,f_1,f_2,f_3,f_4$ must
be 2,2,2,2,2 or 2,2,2,2,3 or 2,2,2,3,3 or 2,2,2,2,4 in some order.
So the common value of the $f_i$ is a sum of two sixth roots of unity.
By rotating by a sixth root of unity, we may assume this value
is 0, 1, $1+\zeta_6$, or 2.
If it is 0 or 1, then the arguments in the proof of
Lemma~\ref{pigeonhole} apply.
Next assume it is $1+\zeta_6$.
The only way two sixth roots of unity can sum to $1+\zeta_6$
is if they are 1 and $\zeta_6$ in some order.
The only ways three sixth roots of unity can sum to $1+\zeta_6$
is if they are $1,1,\zeta_6^2$ or $\zeta_6,\zeta_6,\zeta_6^{-1}$.
So if the numbers of roots of unity occurring in $f_0,f_1,f_2,f_3,f_4$
are 2,2,2,2,2 or 2,2,2,2,3, then $S$ will contain $R_5$ or its rotation
by $\zeta_6$, and the same will be true for 2,2,2,3,3 unless the
two $f_i$ with three terms are $1+1+\zeta_6^2$
and $\zeta_6+\zeta_6+\zeta_6^{-1}$, in which case $S$
contains $(R_5:R_3)$.
It is impossible to write $1+\zeta_6$ as a sum of sixth roots
of unity without using 1 or $\zeta_6$, so if the numbers
are 2,2,2,2,4, then again $S$ contains $R_5$ or its rotation by $\zeta_6$.
Thus we get no new relations where the common value of
the $f_i$ is $1+\zeta_6$.
Lastly, assume this common value is $2$.
Any representation of $2$ as a sum of four or fewer sixth roots
of unity contains $1$, unless
it is $\zeta_6+\zeta_6+\zeta_6^{-1}+\zeta_6^{-1}$, so $S$ will
contain $R_5$ except possibly in the case where $f_0,f_1,f_2,f_3,f_4$
are 2,2,2,2,4 in some order, and the 4 is as above.
But in this final remaining case, $S$ contains $(R_5:R_3)$.
Thus there are no minimal relations $S$ with $p_s=5$ and $10 \le w(S) \le 12$.

\medskip
\noindent{\em Case 3: } $p_s=7$

Since $w(S) \le 12 < 2 \cdot 7$, we can apply Lemma~\ref{pigeonhole}.
Now the sum of $w(T_i)-2$ is required to be $w(S)-7$ which
is at most 5, so the $T$'s that may be used
are $R_3$, $R_5$, $(R_5:R_3)$, and the two of type $(R_5:2R_3)$,
for which weight minus 2 equals 1, 3, 4, and 5, respectively.
So the problem is reduced to listing the partitions of $w(S)-7$
into parts of size 1, 3, 4, and 5.

If all parts used are 1, then we get $(R_7:jR_3)$ with $j=w(S)-7$,
and there are ${7 \choose j}/7$ distinct relations in this class.
Otherwise exactly one part of size 3, 4, or 5 is used, and the
possibilities are as follows.
If a part of size 3 is used, we get $(R_7:R_5)$, $(R_7:R_5,R_3)$,
or $(R_7:R_5,2R_3)$, of weights 10, 11, 12 respectively.
By rotation, the $R_5$ may be assumed to share the 1 in the $R_7$,
and then there are $6 \choose i$ ways to place the $R_3$'s where $i$
is the number of $R_3$'s.
If a part of size 4 is used, we get $(R_7:(R_5:R_3))$ of weight 11
or $(R_7:(R_5:R_3),R_3)$ of weight 12.
By rotation, the $(R_5:R_3)$ may be assumed to share the 1 in the $R_7$,
but any of the six roots of unity in the $(R_5:R_3)$ may be rotated to be 1.
The $R_3$ can then overlap any of the other $6$ seventh roots of unity.
Finally, if a part of size 5 is used, we get $(R_7:(R_5:2R_3))$.
There are two different relations of type $(R_5:2R_3)$ that may be used,
and each has seven roots of unity which may be rotated to be the~1 shared
by the $R_7$, so there are 14 of these all together.

\medskip
\noindent{\em Case 4: } $p_s=11$

Applying Lemma~\ref{pigeonhole} shows that the only possibilities
are $R_{11}$ of weight 11, and $(R_{11}:R_3)$ of weight 12.
\end{pf*}

Now a general relation of weight 12 is a sum of the minimal ones of
weight up to 12, and we can classify them according to the weights
of the minimal relations, which form a partition of 12 with no parts
of size 1 or 4.
We will use the notation $(R_5:2R_3)+2R_3$, for example, to denote
a sum of three minimal relations of type $(R_5:2R_3)$, $R_3$, and $R_3$.
Table~\ref{partitions} lists the possibilities.
The parts may be rotated independently, so any category
involving more than one minimal relation contains
infinitely many relations, even up to rotation (of the entire relation).
Also, the categories are not mutually exclusive, because of the
non-uniqueness of the decomposition into minimal relations.

\begin{table}
\centerline{
\begin{tabular}{c|c}
    Partition	& Relation type		\\ \hline \hline
	12	& $(R_7:5R_3)$		\\ \cline{2-2}
		& $(R_7:R_5,2R_3)$	\\ \cline{2-2}
		& $(R_7:(R_5:R_3),R_3)$	\\ \cline{2-2}
		& $(R_7:(R_5:2R_3))$	\\ \cline{2-2}
		& $(R_{11}:R_3)$	\\ \hline
	10,2	& $(R_7:3R_3)+R_2$	\\ \cline{2-2}
		& $(R_7:R_5)+R_2$	\\ \hline
	9,3	& $(R_5:4R_3)+R_3$	\\ \cline{2-2}
		& $(R_7:2R_3)+R_3$	\\ \hline
	8,2,2	& $(R_5:3R_3)+2R_2$	\\ \cline{2-2}
		& $(R_7:R_3)+2R_2$	\\ \hline
\end{tabular}
\hfil
\begin{tabular}{c|c}
    Partition	& Relation type		\\ \hline \hline
	7,5	& $(R_5:2R_3)+R_5$	\\ \cline{2-2}
		& $R_7+R_5$		\\ \hline
	7,3,2	& $(R_5:2R_3)+R_3+R_2$	\\ \cline{2-2}
		& $R_7+R_3+R_2$		\\ \hline
	6,6	& $2(R_5:R_3)$		\\ \hline
	6,3,3	& $(R_5:R_3)+2R_3$	\\ \hline
	6,2,2,2	& $(R_5:R_3)+3R_2$	\\ \hline
	5,5,2	& $2R_5+R_2$		\\ \hline
	5,3,2,2	& $R_5+R_3+2R_2$	\\ \hline
	3,3,3,3	& $4R_3$		\\ \hline
      3,3,2,2,2	& $2R_3+3R_2$		\\ \hline
    2,2,2,2,2,2	& $6R_2$		\\ \hline
\end{tabular}
}
\vskip12pt
\caption{The types of relations of weight 12.}
\label{partitions}
\end{table}

%****************************************************************************
\section{Solutions to the trigonometric equation}
\label{solutions}

\begin{figure}
\centerline{\psfig{file=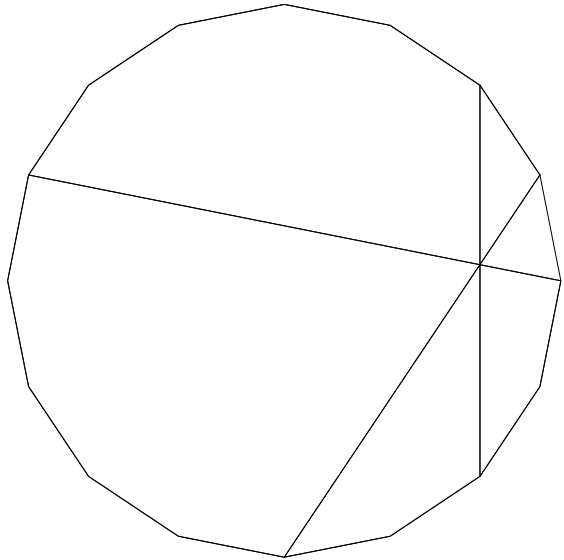,width=3in}}
\caption{A surprising trivial solution for the 16-gon.  The intersection 
point does not lie on any of the 16 lines of symmetry of the
16-gon.\null\hfill\break\null}
\label{16gon}
\end{figure}

Here we use the classification of the previous section to give a complete
listing of the solutions to the trigonometric equation~(\ref{product2}).
There are some obvious solutions to~(\ref{product2}), namely those in
which $U,V,W$ are arbitary positive rational numbers with sum $1/2$,
and $X,Y,Z$ are a permutation of $U,V,W$.
We will call these the trivial solutions, even though the three-diagonal
intersections they give rise to can look surprising.
See Figure~\ref{16gon} for an example on the 16-gon.

The twelve roots of unity occurring in~(\ref{thetwelve}) are not arbitrary;
therefore we must go through Table~\ref{partitions} to see which
relations are of the correct form, i.e., expressible as a sum of six roots
of unity and their inverses, where the product of the six is -1.
First let us prove a few lemmas that will greatly reduce the number of cases.

\begin{lemma}
Let $S$ be a relation of weight $k \le 12$.
Suppose $S$ is stable under complex conjugation (i.e., under
$\zeta \mapsto \zeta^{-1}$).
Then $S$ has a complex conjugation-stable decomposition
into minimal relations; i.e., each minimal relation occurring is itself
stable under complex conjugation, or can be paired with another minimal
relation which is its complex conjugate.
\end{lemma}

\begin{pf}
We will use induction on $k$.
If $S$ is minimal, there is nothing to prove.
Otherwise let $T$ be a (minimal) subrelation of $S$ of minimal weight,
so $T$ is of weight at most 6.
The complex conjugate $\overline{T}$ of $T$ is another minimal relation
in $S$.
If they do not intersect, then we take the decomposition of $S$ into $T$,
$\overline{T}$, and a decomposition of $S \setminus (T \cup \overline{T})$
given by the inductive hypothesis.
If they do overlap and the weight of $T$ is at most 5, then $T=R_p$ for
some prime $p$, and the fact that $T$ intersects $\overline{T}$ implies
that $T=\overline{T}$, and we get the result by applying the inductive
hypothesis to $S \setminus T$.

The only remaining case is where $S$ is of type $2(R_5:R_3)$.
If the two $(R_5:R_3)$'s are not conjugate to each other, then for each
there is a root of unity $\zeta$ such that $\zeta$ and $\zeta^{-1}$ occur
in that (rotation of) $(R_5:R_3)$.
The quotient $\zeta^2$ is then a $30$-th root of unity, so $\zeta$
itself is a $60$-th root of unity.
Thus each $(R_5:R_3)$ is a rotation of the ``standard'' $(R_5:R_3)$
as in~(\ref{r5r3}) by a $60$-th root of unity, and we let Mathematica
check the $60^2$ possibilities.
\end{pf}

We do not know if the preceding lemma holds for relations of weight
greater than 12.

\begin{lemma}
Let $S$ be a minimal relation of type $(R_p:T_1,\ldots,T_j)$, $p \ge 5$,
where the $T_i$ involve roots of unity of order prime to $p$, and $j<p$.
If $S$ is stable under complex conjugation, then the particular rotation
of $R_p$ from which the $T_i$ were ``subtracted'' is also stable
(and hence so is the collection of the relations subtracted).
\end{lemma}

\begin{pf}
Let $\ell$ be the product of the orders of the roots of unity in
all the $T_i$.
The elements of $S$ in the original $R_p$ can be characterized
as those terms of $S$ that are unique in their coset of $\mu_\ell$
(the $\ell$-th roots of unity), and this condition is stable under
complex conjugation, so the set of terms of the $R_p$ that were not
subtracted is stable.
Since $j<p$, we can pick one such term $\zeta$.
Then the quotient $\zeta/\zeta^{-1}$ is a $p$-th root of unity,
so $\zeta$ is a $2p$-th root of unity, and hence the $R_p$
containing it is stable.
\end{pf}

\begin{cor}
A relation of type $(R_7:(R_5:R_3),R_3)$ cannot be stable under
complex conjugation.
\end{cor}

Even with these restrictions, a very large number of cases remain,
so we perform the calculation using Mathematica.
Each entry of Table~\ref{partitions} represents a finite number of
linearly parameterized (in the exponents) families of relations of weight 12.
For each parameterized family, we check to see what additional
constraints must be put on the parameters for the relation to be
of the form of~(\ref{thetwelve}).
Next, for each parameterized family of solutions to~(\ref{thetwelve}),
we calculate the corresponding $U,V,W,X,Y,Z$ and throw away solutions
in which some of these are nonpositive.
Finally, we sort $U,V,W$ and $X,Y,Z$ and interchange the two triples
if $U>X$, in order to count the solutions only up to symmetry.

The results of this computation are recorded in the following theorem.

\begin{theorem}
\label{thesolutions}
The positive rational solutions to~(\ref{product2}), up to symmetry,
can be classified as follows:
\begin{enumerate}
	\item The trivial solutions,
		which arise from relations of type $6R_2$.
	\item Four one-parameter families of solutions,
		listed in Table~\ref{families}.
		The first arises from relations of type $4R_3$,
		and the other three arise from relations of type $2R_3+3R_2$.
	\item Sixty-five ``sporadic'' solutions, listed in
		Table~\ref{sporadics}, which arise from the other types
		of weight 12 relations listed in Table~\ref{partitions}.
\end{enumerate}
The only duplications in this list are that the second family of
Table~\ref{families} gives a trivial solution for $t=1/12$,
the first and fourth families of Table~\ref{families} give the same
solution when $t=1/18$ in both, and the second and fourth families
of Table~\ref{families} give the same solution when $t=1/24$ in both.
\end{theorem}

\begin{table}
\begin{center}
\begin{tabular}{|c|c|c||c|c|c||c|}
 $U$ & $V$ & $W$ & $X$ & $Y$ & $Z$ & Range \\ \hline \hline
 $1/6$ & $t$ & $1/3-2t$ & $1/3+t$ & $t$ & $1/6-t$ & $0<t<1/6$ \\ \hline
 $1/6$ & $1/2-3t$ & $t$ & $1/6-t$ & $2t$ & $1/6+t$ & $0<t<1/6$ \\ \hline
 $1/6$ & $1/6-2t$ & $2t$ & $1/6-2t$ & $t$ & $1/2+t$ & $0<t<1/12$ \\ \hline
 $1/3-4t$ & $t$ & $1/3+t$ & $1/6-2t$ & $3t$ & $1/6+t$ & $0<t<1/12$ \\ \hline
\end{tabular}
\end{center}
\caption{The nontrivial infinite families of solutions to~(\ref{product2}).}
\label{families}
\end{table}

\begin{table}
{\tiny
\begin{center}
\begin{tabular}{c||c|c|c||c|c|c||c}
Denominator	& $U$ & $V$ & $W$ & $X$ & $Y$ & $Z$
		& Relation type \\ \hline \hline
30 & $1/10$ & $2/15$ & $3/10$ & $2/15$ & $1/6$ & $1/6$
		& $2(R_5:R_3)$ \\ \cline{2-7}
 & $1/15$ & $1/15$ & $7/15$ & $1/15$ & $1/10$ & $7/30$ & \\ \cline{2-7}
 & $1/30$ & $7/30$ & $4/15$ & $1/15$ & $1/10$ & $3/10$ & \\ \cline{2-7}
 & $1/30$ & $1/10$ & $7/15$ & $1/15$ & $1/15$ & $4/15$ & \\ \cline{2-7}
 & $1/30$ & $1/15$ & $19/30$ & $1/15$ & $1/10$ & $1/10$ & \\ \cline{2-8}
 & $1/15$ & $1/6$ & $4/15$ & $1/10$ & $1/10$ & $3/10$
		& $(R_5:R_3)+2R_3$ \\ \cline{2-7}
 & $1/15$ & $2/15$ & $11/30$ & $1/10$ & $1/6$ & $1/6$ & \\ \cline{2-7}
 & $1/30$ & $1/6$ & $13/30$ & $1/10$ & $2/15$ & $2/15$ & \\ \cline{2-7}
 & $1/30$ & $1/30$ & $7/10$ & $1/30$ & $1/15$ & $2/15$ & \\ \cline{2-8}
 & $1/30$ & $7/30$ & $3/10$ & $1/15$ & $2/15$ & $7/30$
		& $R_5+R_3+2R_2$ \\ \cline{2-7}
 & $1/30$ & $1/6$ & $11/30$ & $1/15$ & $1/10$ & $4/15$ & \\ \cline{2-7}
 & $1/30$ & $1/10$ & $13/30$ & $1/30$ & $2/15$ & $4/15$ & \\ \cline{2-7}
 & $1/30$ & $1/15$ & $8/15$ & $1/30$ & $1/10$ & $7/30$ & \\ \hline
42 & $1/14$ & $5/42$ & $5/14$ & $2/21$ & $5/42$ & $5/21$
		& $(R_7:5R_3)$ \\ \cline{2-7}
 & $1/21$ & $4/21$ & $13/42$ & $1/14$ & $1/6$ & $3/14$ & \\ \cline{2-7}
 & $1/42$ & $3/14$ & $5/14$ & $1/21$ & $1/6$ & $4/21$ & \\ \cline{2-7}
 & $1/42$ & $1/6$ & $19/42$ & $1/14$ & $2/21$ & $4/21$ & \\ \cline{2-7}
 & $1/42$ & $1/6$ & $13/42$ & $1/21$ & $1/14$ & $8/21$ & \\ \cline{2-7}
 & $1/42$ & $1/21$ & $13/21$ & $1/42$ & $1/14$ & $3/14$ & \\ \hline
60 & $1/20$ & $1/12$ & $29/60$ & $1/15$ & $1/10$ & $13/60$
		& $2(R_5:R_3)$ \\ \cline{2-7}
 & $1/20$ & $1/12$ & $9/20$ & $1/15$ & $1/12$ & $4/15$ & \\ \cline{2-7}
 & $1/20$ & $1/12$ & $5/12$ & $1/20$ & $1/10$ & $3/10$ & \\ \cline{2-7}
 & $1/60$ & $4/15$ & $3/10$ & $1/20$ & $1/12$ & $17/60$ & \\ \cline{2-7}
 & $1/60$ & $13/60$ & $9/20$ & $1/12$ & $1/10$ & $2/15$ & \\ \cline{2-7}
 & $1/60$ & $13/60$ & $5/12$ & $1/20$ & $2/15$ & $1/6$ & \\ \cline{2-8}
 & $1/12$ & $1/6$ & $17/60$ & $2/15$ & $3/20$ & $11/60$
		& $(R_5:3R_3)+2R_2$ \\ \cline{2-7}
 & $1/12$ & $2/15$ & $19/60$ & $1/10$ & $3/20$ & $13/60$ & \\ \cline{2-7}
 & $1/15$ & $11/60$ & $13/60$ & $1/12$ & $1/10$ & $7/20$ & \\ \cline{2-7}
 & $1/20$ & $11/60$ & $3/10$ & $1/12$ & $7/60$ & $4/15$ & \\ \cline{2-7}
 & $1/20$ & $1/10$ & $23/60$ & $1/15$ & $1/12$ & $19/60$ & \\ \cline{2-7}
 & $1/30$ & $7/60$ & $19/60$ & $1/20$ & $1/15$ & $5/12$ & \\ \cline{2-7}
 & $1/30$ & $1/12$ & $7/12$ & $1/15$ & $1/10$ & $2/15$ & \\ \cline{2-7}
 & $1/30$ & $1/20$ & $11/20$ & $1/30$ & $1/15$ & $4/15$ & \\ \cline{2-7}
 & $1/60$ & $3/10$ & $7/20$ & $1/12$ & $7/60$ & $2/15$ & \\ \cline{2-7}
 & $1/60$ & $4/15$ & $23/60$ & $1/12$ & $1/10$ & $3/20$ & \\ \cline{2-7}
 & $1/60$ & $7/30$ & $5/12$ & $1/15$ & $7/60$ & $3/20$ & \\ \cline{2-7}
 & $1/60$ & $13/60$ & $11/30$ & $1/20$ & $1/12$ & $4/15$ & \\ \cline{2-7}
 & $1/60$ & $1/6$ & $31/60$ & $1/15$ & $1/10$ & $2/15$ & \\ \cline{2-7}
 & $1/60$ & $1/6$ & $5/12$ & $1/20$ & $1/15$ & $17/60$ & \\ \cline{2-7}
 & $1/60$ & $2/15$ & $9/20$ & $1/30$ & $1/12$ & $17/60$ & \\ \cline{2-7}
 & $1/60$ & $1/10$ & $31/60$ & $1/30$ & $1/15$ & $4/15$ & \\ \hline
84 & $1/12$ & $3/14$ & $19/84$ & $11/84$ & $13/84$ & $4/21$
		& $(R_7:R_3)+2R_2$ \\ \cline{2-7}
 & $1/14$ & $11/84$ & $23/84$ & $1/12$ & $2/21$ & $29/84$ & \\ \cline{2-7}
 & $1/21$ & $13/84$ & $23/84$ & $1/14$ & $1/12$ & $31/84$ & \\ \cline{2-7}
 & $1/42$ & $1/12$ & $7/12$ & $1/21$ & $1/14$ & $4/21$ & \\ \cline{2-7}
 & $1/84$ & $25/84$ & $5/14$ & $5/84$ & $1/12$ & $4/21$ & \\ \cline{2-7}
 & $1/84$ & $5/21$ & $5/12$ & $5/84$ & $1/14$ & $17/84$ & \\ \cline{2-7}
 & $1/84$ & $3/14$ & $37/84$ & $1/21$ & $1/12$ & $17/84$ & \\ \cline{2-7}
 & $1/84$ & $1/6$ & $43/84$ & $1/21$ & $1/14$ & $4/21$ & \\ \hline
90 & $1/18$ & $13/90$ & $7/18$ & $11/90$ & $2/15$ & $7/45$
		& $(R_5:R_3)+2R_3$ \\ \cline{2-7}
 & $1/45$ & $19/90$ & $16/45$ & $1/18$ & $1/10$ & $23/90$ & \\ \cline{2-7}
 & $1/90$ & $23/90$ & $31/90$ & $2/45$ & $1/15$ & $5/18$ & \\ \cline{2-7}
 & $1/90$ & $17/90$ & $47/90$ & $1/18$ & $4/45$ & $2/15$ & \\ \hline
120 & $13/120$ & $3/20$ & $31/120$ & $2/15$ & $19/120$ & $23/120$ &
		$(R_5:R_3)+3R_2$ \\ \cline{2-7}
 & $1/12$ & $19/120$ & $29/120$ & $1/10$ & $13/120$ & $37/120$ & \\ \cline{2-7}
 & $1/20$ & $23/120$ & $29/120$ & $1/15$ & $13/120$ & $41/120$ & \\ \cline{2-7}
 & $1/60$ & $13/120$ & $73/120$ & $1/20$ & $1/12$ & $2/15$ & \\ \cline{2-7}
 & $1/120$ & $7/20$ & $43/120$ & $7/120$ & $11/120$ & $2/15$ & \\ \cline{2-7}
 & $1/120$ & $3/10$ & $49/120$ & $7/120$ & $1/12$ & $17/120$ & \\ \cline{2-7}
 & $1/120$ & $4/15$ & $53/120$ & $1/20$ & $11/120$ & $17/120$ & \\ \cline{2-7}
 & $1/120$ & $13/60$ & $61/120$ & $1/20$ & $1/12$ & $2/15$ & \\ \hline
210 & $1/15$ & $41/210$ & $8/35$ & $1/14$ & $31/210$ & $61/210$
		& $(R_7:(R_5:2R_3))$ \\ \cline{2-7}
 & $13/210$ & $1/10$ & $83/210$ & $1/14$ & $4/35$ & $9/35$ & \\ \cline{2-7}
 & $1/35$ & $2/15$ & $97/210$ & $1/14$ & $17/210$ & $47/210$ & \\ \cline{2-7}
 & $1/210$ & $3/14$ & $121/210$ & $11/210$ & $1/15$ & $3/35$ & \\ \hline
\end{tabular}
\end{center}
}
\caption{The 65 sporadic solutions to~(\ref{product2}).}
\label{sporadics}
\end{table}

Some explanation of the tables is in order.
The last column of Table~\ref{families} gives the allowable range
for the rational parameter $t$.
The entries of Table~\ref{sporadics} are sorted according to the least
common denominator of $U,V,W,X,Y,Z$, which is also the least $n$
for which diagonals of a regular $n$-gon can create arcs of the
corresponding lengths.
The relation type from which each solution derives is also given.
The reason 11 does not appear in the least common denominator for
any sporadic solution is that the relation $(R_{11}:R_3)$
cannot be put in the form of~(\ref{thetwelve}) with the $\alpha_j$
summing to 1, and hence leads to no solutions of~(\ref{product2}).
(Several other types of relations also give rise to no solutions.)

Tables~\ref{families} and~\ref{sporadics} are the same as Bol's tables
at the bottom of page~40 and on page~41 of~\cite{bol},
in a slightly different format.

\medskip

The arcs cut by diagonals of a regular $n$-gon have lengths which are
multiples of $2 \pi/n$, so $U$, $V$, $W$, $X$, $Y$ and $Z$ corresponding
to any configuration of three diagonals meeting must be multiples of $1/n$.
With this additional restriction, trivial solutions to~(\ref{product2})
occur only when $n$ is even (and at least 6).
Solutions within the infinite families of Table~\ref{families} occur
when $n$ is a multiple of 6 (and at least 12), and there $t$ must be
a multiple of $1/n$.
Sporadic solutions with least common denominator $d$ occur if and only
if $n$ is a multiple of $d$.

%****************************************************************************
\section{Intersections of more than three diagonals}
\label{more}

Now that we know the configurations of three diagonals meeting,
we can check how they overlap to produce configurations of more
than three diagonals meeting.
We will disregard configurations in which the intersection point
is the center of the $n$-gon, since these are easily described:
there are exactly $n/2$ diagonals (diameters) through the center
when $n$ is even, and none otherwise.

When $k$ diagonals meet, they form $2k$ arcs, whose lengths we will
measure as a fraction of the whole circumference (so they will be
multiples of $1/n$) and list in counterclockwise order.
(Warning: this is different from the order used in Tables \ref{families}
and~\ref{sporadics}.)
The least common denominator of the numbers in this list will be
called the denominator of the configuration.
It is the least $n$ for which the configuration can be realized as
diagonals of a regular $n$-gon.

\begin{lemma}
\label{denominator}
If a configuration of $k \ge 2$ diagonals meeting at an interior point
other than the center has denominator dividing $d$, then any configuration
of diagonals meeting at that point has denominator dividing $\LCM(2d,3)$.
\end{lemma}

\begin{pf}
We may assume $k=2$.
Any other configuration of diagonals through the intersection point is
contained in the union of configurations obtained by adding one
diagonal to the original two, so we may assume the final configuration
consists of three diagonals, two of which were the original two.
Now we need only go through our list of three-diagonal intersections.

It can be checked (using Mathematica) that removing any diagonal from
a sporadic configuration of three intersecting diagonals yields a
configuration whose denominator is the same or half as much, except
that it is possible that removing a diagonal from a three-diagonal
configuration of denominator 210 or 60 yields one of denominator 70 or 20,
respectively, which proves the desired result for these cases.
The additive group generated by $1/6$ and the normalized arc lengths
of a configuration obtained by removing a diagonal from a configuration
corresponding to one of the families of Table~\ref{families}
contains $2t$ where $t$ is the parameter, (as can be verified using
Mathematica again), which means that adding that third diagonal
can at most double the denominator (and throw in a factor of 3,
if it isn't already there).
Similarly, it is easily checked (even by hand), that the subgroup
generated by the normalized arc lengths of a configuration obtained
by removing one of the three diagonals of a configuration corresponding
to a trivial solution to~(\ref{product2}) but with intersection point
not the center, contains twice the arc lengths of the original configuration.
\end{pf}

\begin{cor}
\label{partsporadic}
If a configuration of three or more diagonals meeting includes three
forming a sporadic configuration, then its
denominator is 30, 42, 60, 84, 90, 120, 168, 180, 210, 240, or 420.
\end{cor}

\begin{pf}
Combine the lemma with the list of denominators of sporadic configurations
listed in Table~\ref{sporadics}.
\end{pf}

For $k \ge 4$, a list of $2k$ positive rational numbers summing to~1
arises this way if and only if the lists of length $2k-2$ which would
arise by removing the first or second diagonal actually
correspond to $k-1$ intersecting diagonals.
Suppose $k=4$.
If we specify the sporadic configuration or parameterized family
of configurations that arise when we remove the first or second diagonal,
we get a set of linear conditions on the eight arc lengths.
Corollary~\ref{partsporadic} tells us that we get a configuration with
denominator among 30, 42, 60, 84, 90, 120, 168, 180, 210, 240, and 420,
if one of these two is sporadic.
Using Mathematica to perform this computation for the rest of
possibilities in Theorem~\ref{thesolutions} shows that the other
four-diagonal configurations, up to rotation and reflection,
fall into 12 one-parameter families, which are listed in
Table~\ref{theeights} by the eight normalized arc lengths
and the range for the parameter $t$, with a finite number of exceptions
of denominators among 12, 18, 24, 30, 36, 42, 48, 60, 84, and 120.

\begin{table}
\centerline{\small
\begin{tabular}{|c|c|c|c|c|c|c|c||c|}
  & & & & & & & & Range \\ \hline \hline
 $t$ & $t$ & $t$ & $1/6-2t$
	& $1/6$ & $1/3+t$ & $1/6$ & $1/6-2t$ & $0<t<1/12$ \\ \hline
 $t$ & $1/6-t$ & $1/6-t$ & $1/6-t$
	& $t$ & $1/6$ & $1/6+t$ & $1/6$ & $0<t<1/6$ \\ \hline
 $1/6-4t$ & $2t$ & $t$ & $3t$
	& $1/6-4t$ & $1/6$ & $1/6+t$ & $1/3+t$ & $0<t<1/24$ \\ \hline
 $2t$ & $1/2-t$ & $2t$ & $1/6-2t$
	& $t$ & $1/6-t$ & $t$ & $1/6-2t$ & $0<t<1/12$ \\ \hline
 $1/3-4t$ & $1/6+t$ & $1/2-3t$ & $-1/6+4t$
	& $1/6-2t$ & $t$ & $1/6-t$ & $-1/6+4t$ & $1/24<t<1/12$ \\ \hline
 $2t$ & $t$ & $3t$ & $1/6-2t$
	& $1/6$ & $1/6-t$ & $1/3-t$ & $1/6-2t$ & $0<t<1/12$ \\ \hline
 $t$ & $t$ & $2t$ & $1/3-t$
	& $1/6$ & $1/6-t$ & $1/6-t$ & $1/6-t$ & $0<t<1/6$ \\ \hline
 $1/3-4t$ & $1/6$ & $t$ & $t$
	& $1/6-2t$ & $1/3-2t$ & $3t$ & $3t$ & $0<t<1/12$ \\ \hline
 $2t$ & $1/3-2t$ & $1/6-t$ & $1/6-t$
	& $1/6$ & $1/6$ & $t$ & $t$ & $0<t<1/6$ \\ \hline
 $1/3-4t$ & $2t$ & $t$ & $t$
	& $1/6-2t$ & $1/6$ & $1/6+t$ & $1/6+t$ & $0<t<1/12$ \\ \hline
 $1/3-4t$ & $2t$ & $1/6-t$ & $t$
	& $1/6-2t$ & $2t$ & $1/3-t$ & $3t$ & $0<t<1/12$ \\ \hline
 $2t$ & $1/6-t$ & $t$ & $1/6-t$
	& $t$ & $1/6-t$ & $2t$ & $1/2-3t$ & $0<t<1/6$ \\ \hline
\end{tabular}
}
\vskip10pt
\caption{The one-parameter families of four-diagonal configurations.}
\label{theeights}
\end{table}

We will use a similar argument when $k=5$.
Any five-diagonal configuration containing a sporadic three-diagonal
configuration will again have denominator
among 30, 42, 60, 84, 90, 120, 168, 180, 210, 240, and 420.
Any other five-diagonal configuration containing one of the exceptional
four-diagonal configurations will have denominator
among 12, 18, 24, 30, 36, 42, 48, 60, 72, 84, 96, 120, 168, and 240,
by Lemma~\ref{denominator}.
Finally, another Mathematica computation shows that the one-parameter
families of four-diagonal configurations overlap to produce the
one-parameter families listed (up to rotation and reflection)
in Table~\ref{thetens}, and a finite number of exceptions of
denominators among 18, 24, and 30.

\begin{table}
\centerline{\small
\begin{tabular}{|c|c|c|c|c|c|c|c|c|c||c|}
  & & & & & & & & & & Range \\ \hline \hline
 $t$ & $2t$ & $1/6-2t$ & $1/6$ & $1/6-t$
	& $1/6-t$ & $1/6$ & $1/6-2t$ & $2t$ & $t$ & $0<t<1/12$ \\ \hline
 $t$ & $2t$ & $1/6-4t$ & $1/6$ & $1/6+t$
	& $1/6+t$ & $1/6$ & $1/6-4t$ & $2t$ & $t$ & $0<t<1/24$ \\ \hline
 $t$ & $1/6-2t$ & $-1/6+4t$ & $1/3-4t$ & $1/6+t$ & $1/6+t$
	& $1/3-4t$ & $-1/6+4t$ & $1/6-2t$ & $t$ & $1/24<t<1/12$ \\ \hline
 $t$ & $1/6-2t$ & $2t$ & $1/3-4t$ & $3t$ & $3t$
	& $1/3-4t$ & $2t$ & $1/6-2t$ & $t$ & $0<t<1/12$ \\ \hline
\end{tabular}
}
\vskip10pt
\caption{The one-parameter families of five-diagonal configurations.}
\label{thetens}
\end{table}

For $k=6$, any six-diagonal configuration containing a sporadic
three-diagonal configuration will again have denominator
among 30, 42, 60, 84, 90, 120, 168, 180, 210, 240, and 420.
Any six-diagonal configuration containing one of the exceptional
four-diagonal configurations will have denominator
among 12, 18, 24, 30, 36, 42, 48, 60, 72, 84, 96, 120, 168, and 240.
Any six-diagonal configuration containing one of the
exceptional five-diagonal configurations will have denominator
among 18, 24, 30, 36, 48, and 60.
Another Mathematica computation shows that the one-parameter families
of five-diagonal configurations cannot combine to give a
six-diagonal configuration.

Finally for $k \ge 7$, any $k$-diagonal configuration must contain
an exceptional configuration of 3, 4, or 5 diagonals, and hence by
Lemma~\ref{denominator} has denominator
among 12, 18, 24, 30, 36, 42, 48, 60, 72, 84, 90, 96,
120, 168, 180, 210, 240, and 420.

We summarize the results of this section in the following.

\begin{prop}
\label{configurations}
The configurations of $k \ge 4$ diagonals meeting at a point
not the center, up to rotation and reflection, fall into the
one-parameter families listed in Tables \ref{theeights} and~\ref{thetens},
with finitely many exceptions (for fixed $k$) of denominators
among 12, 18, 24, 30, 36, 42, 48, 60, 72, 84, 90, 96, 120, 168,
180, 210, 240, and 420.
\end{prop}
In fact, many of the numbers listed in the proposition do not
actually occur as denominators of exceptional configurations.
For example, it will turn out that the only denominator greater
than 120 that occurs is 210.

%****************************************************************************
\section{The formula for intersection points}
\label{points}

Let $a_k(n)$ denote the number of points inside the regular $n$-gon
other than the center where exactly $k$ lines meet.
Let $b_k(n)$ denote the number of $k$-tuples of diagonals which meet
at a point inside the $n$-gon other than the center.
Each interior point at which exactly $m$ diagonals meet gives rise
to $m \choose k$ such $k$-tuples, so we have the relationship
\begin{equation}
\label{relationship}
	b_k(n) = \sum_{m \ge k} {m \choose k} a_m(n)
\end{equation}
Since every four distinct vertices of the $n$-gon determine one pair
of diagonals which intersect inside, the number of such pairs is
exactly $n \choose 4$, but if $n$ is even, then ${n/2} \choose 2$
of these are pairs which meet at the center, so
\begin{equation}
\label{b2}
	b_2(n)={n \choose 4} - {{n/2} \choose 2} \delta_2(n).
\end{equation}
(Recall that $\delta_m(n)$ is defined to be 1 if $n$ is a multiple of $m$,
and 0 otherwise.)

We will use the results of the previous two sections to deduce the
form of $b_k(n)$ and then the form of $a_k(n)$.
To avoid having to repeat the following, let us make a definition.

\begin{defn}
A function on integers $n \ge 3$ will be called {\em tame}
if it is a linear combination (with rational coefficients)
of the functions $n^3$, $n^2$, $n$, $1$, $n^2 \delta_2(n)$,
$n \delta_2(n)$, $\delta_2(n)$, $\delta_4(n)$, $n \delta_6(n)$,
$\delta_6(n)$, $\delta_{12}(n)$, $\delta_{18}(n)$, $\delta_{24}(n)$,
$\delta_{24}(n-6)$, $\delta_{30}(n)$, $\delta_{36}(n)$, $\delta_{42}(n)$,
$\delta_{48}(n)$, $\delta_{60}(n)$, $\delta_{72}(n)$, $\delta_{84}(n)$,
$\delta_{90}(n)$, $\delta_{96}(n)$, $\delta_{120}(n)$, $\delta_{168}(n)$,
$\delta_{180}(n)$, $\delta_{210}(n)$, and $\delta_{420}(n)$.
\end{defn}

\begin{prop}
\label{form}
For each $k \ge 2$, the function $b_k(n)/n$ on integers $n \ge 3$ is tame.
\end{prop}

\begin{pf}
The case $k=2$ is handled by~(\ref{b2}), so assume $k \ge 3$.
Each list of $2k$ normalized arc lengths as in Section~\ref{more}
corresponding to a configuration of $k$ diagonals meeting at a point
other than the center, considered up to rotation (but not reflection),
contributes $n$ to $b_k(n)$.
(There are $n$ places to start measuring the arcs from, and these $n$
configurations are distinct, because the corresponding intersection points
differ by rotations of multiples of $2 \pi/n$, and by assumption they are
not at the center.)
So $b_k(n)/n$ counts such lists.

Suppose $k=3$.
When $n$ is even, the family of trivial solutions to the trigonometric
equation~(\ref{product2}) has $U=a/n$, $V=b/n$, $W=c/n$,
where $a$, $b$, and $c$ are positive integers with sum $n/2$,
and $X$, $Y$, and $Z$ are some permutation of $U$, $V$, $W$.
Each permutation gives rise to a two-parameter family of six-long lists
of arc lengths, and the number of lists within each family is the number
of partitions of $n/2$ into three positive parts, which is a quadratic
polynomial in $n$.
Similarly each family of solutions in Table~\ref{families} gives rise
to a number of one-parameter families of lists, when $n$ is a multiple of 6,
each containing $\lceil n/6 \rceil - 1$ or $\lceil n/12 \rceil - 1$ lists.
These functions of $n$ (extended to be 0 when 6 does not divide $n$)
are expressible as a linear combination of $n \delta_6(n)$, $\delta_6(n)$,
and $\delta_{12}(n)$.
Finally the sporadic solutions to~\ref{product2} give rise to a
finite number of lists, having denominators among 30, 42, 60, 84, 90, 120,
and 210, so their contribution to $b_3(n)/n$ is a linear combination
of $\delta_{30}(n), \ldots, \delta_{210}(n)$.

But these families of lists overlap, so we must use the Principle of
Inclusion-Exclusion to count them properly.
To show that the result is a tame function, it suffices to show that
the number of lists in any intersection of these families is a tame function.
When two of the trivial families overlap but do not coincide, they
overlap where two of the $a$, $b$, and $c$ above are equal, and the
corresponding lists lie in one of the one-parameter families
$(t,t,t,t,1/2-2t,1/2-2t)$ or $(t,t,t,1/2-2t,t,1/2-2t)$ (with $0<t<1/4$),
each of which contain $\lceil n/4 \rceil -1$ lists (for $n$ even).
This function of $n$ is a combination of $n \delta_2(n)$, $\delta_2(n)$,
and $\delta_4(n)$, hence it is tame.
Any other intersection of the infinite families must contain the
intersection of two one-parameter families which are among the two
above or arise from Table~\ref{families}, and a Mathematica
computation shows that such an intersection consists of at most a
single list of denominator among 6, 12, 18, 24, and 30.
And, of course, any intersection involving a single sporadic list,
can contain at most that sporadic list.
Thus the number of lists within any intersection is a tame function
of $n$.
Finally we must delete the lists which correspond to configurations
of diagonals meeting at the center.
These are the lists within the trivial two-parameter family
$(t,u,1/2-t-u,t,u,1/2-t-u)$, so their number is also a tame function
of $n$, by the Principle of Inclusion-Exclusion again.
Thus $b_3(n)/n$ is tame.

Next suppose $k=4$.
The number of lists within each family listed in Table~\ref{theeights},
or the reflection of such a family, is (when $n$ is divisible by 6)
the number of multiples of $1/n$ strictly between $\alpha$ and $\beta$,
where the range for the parameter $t$ is $\alpha<t<\beta$.
This number is $\lceil \beta n \rceil - 1 - \lfloor \alpha n \rfloor$.
Since the table shows that  $\alpha$ and $\beta$ are always multiples
of $1/24$, this function of $n$ is expressible as a combination
of $n \delta_6(n)$ and a function on multiples of 6 depending only
on $n \bmod 24$, and the latter can be written as a combination
of $\delta_6(n)$, $\delta_{12}(n)$, $\delta_{24}(n)$,
and $\delta_{24}(n-6)$, so it is tame.
Mathematica shows that when two of these families are not the same,
they intersect in at most a single list of denominator
among 6, 12, 18, and 24.
So these and the exceptions of Proposition~\ref{configurations}
can be counted by a tame function.
Thus, again by the Principle of Inclusion-Exclusion, $b_4(n)/n$ is tame.

The proof for $k=5$ is identical to that of $k=4$, using
Table~\ref{thetens} instead of Table~\ref{theeights},
and using another Mathematica computation which shows that
the intersections of two one-parameter families of lists consist
of at most a single list of denominator 24.

The proof for $k \ge 6$ is even simpler, because then there are
only the exceptional lists.
By Proposition~\ref{configurations}, $b_k(n)/n$ is a linear
combination of $\delta_m(n)$ where $m$ ranges over the possible
denominators of exceptional lists listed in the proposition, so it is tame.
\end{pf}

\begin{lemma}
\label{determined}
A tame function is determined by its values
at $n=$ 3, 4, 5, 6, 7, 8, 9, 10, 12, 18, 24, 30, 36, 42, 48, 54, 60,
66, 72, 84, 90, 96, 120, 168, 180, 210, and 420.
\end{lemma}

\begin{pf}
By linearity, it suffices to show that if a tame function $f$
is zero at those values, then $f$ is the zero linear combination
of the functions in the definition of a tame function.
The vanishing at $n=3$, 5, 7, and 9 forces the coefficients
of $n^3$, $n^2$, $n$, and $1$ to vanish, by Lagrange interpolation.
Then comparing the values at $n=4$ and $n=10$ shows that the
coefficient of $\delta_4(n)$ is zero.
The vanishing at $n=4$, 8, and 10 forces the coefficients
of $n^2 \delta_2(n)$, $n \delta_2(n)$, and $\delta_2(n)$ to vanish.
Comparing the values at $n=6$ and $n=54$ shows that the
coefficient of $n \delta_6$ is zero.
Comparing the values at $n=6$ and $n=66$ shows that the
coefficient of $\delta_24(n-6)$ is zero.

At this point, we know that $f(n)$ is a combination of $\delta_m(n)$,
for $m=6$, 12, 18, 24, 30, 36, 42, 48, 60, 72, 84, 90, 96, 120,
168, 180, 210, and 420.
For each $m$ in turn, $f(m)=0$ now implies that the coefficient
of $\delta_m(n)$ is zero.
\end{pf}

\begin{pf*}{Proof of Theorem~\ref{countintersections}}
Computation (see the appendix) shows that the tame function $b_8(n)/n$
vanishes at all the numbers listed in Lemma~\ref{determined}.
Hence by that lemma, $b_8(n)=0$ for all $n$.
Thus by~(\ref{relationship}), $a_k(n)$ and $b_k(n)$ are identically zero
for all $k \ge 8$ as well.

By reverse induction on $k$, we can invert~(\ref{relationship})
to express $a_k(n)$ as a linear combination of $b_m(n)$ with $m \ge k$.
Hence $a_k(n)/n$ is tame as well for each $k \ge 2$.
Computation shows that the equations
\begin{eqnarray*}
	a_2(n)/n	& = &	(n^3 - 6n^2 + 11n - 6)/24
		+ (-5 n^2 + 46 n - 72)/16 \cdot \delta_2(n) \\
			&   &	\mbox{} - 9/4 \cdot \delta_4(n)
		+ (-19 n + 110)/2 \cdot \delta_6(n)
		+ 54 \cdot \delta_{12}(n) + 84 \cdot \delta_{18}(n)	\\
			&   &	\mbox{} + 50 \cdot \delta_{24}(n)
		- 24 \cdot \delta_{30}(n) - 100 \cdot \delta_{42}(n)
		- 432 \cdot \delta_{60}(n)	\\
			&   &	\mbox{} -204 \cdot \delta_{84}(n)
		- 144 \cdot \delta_{90}(n) - 204 \cdot \delta_{120}(n)
		- 144 \cdot \delta_{210}(n)	\\
	a_3(n)/n	& = &	(5 n^2 - 48 n + 76)/48 \cdot \delta_2(n)
		+ 3/4 \cdot \delta_4(n) + (7n - 38)/6 \cdot \delta_6(n)\\
			&   &	\mbox{} - 8 \cdot \delta_{12}(n)
		- 20 \cdot \delta_{18}(n) - 16 \cdot \delta_{24}(n)
		- 19 \cdot \delta_{30}(n) + 8 \cdot \delta_{42}(n)	\\
			&   &	\mbox{} + 68 \cdot \delta_{60}(n)
		+ 60 \cdot \delta_{84}(n) + 48 \cdot \delta_{90}(n)
		+ 60 \cdot \delta_{120}(n) + 48 \cdot \delta_{210}(n)	\\
	a_4(n)/n	& = &	(7n - 42)/12 \cdot \delta_6(n)
		- 5/2 \cdot \delta_{12}(n) - 4 \cdot \delta_{18}(n)
		+ 3 \cdot \delta_{24}(n)	\\
			&   &	\mbox{} + 6 \cdot \delta_{42}(n)
		+ 34 \cdot \delta_{60}(n) - 6 \cdot \delta_{84}(n)
		- 6 \cdot \delta_{120}(n)	\\
	a_5(n)/n	& = &	(n - 6)/4 \cdot \delta_6(n)
		- 3/2 \cdot \delta_{12}(n) - 2 \cdot \delta_{24}(n)
		+ 4 \cdot \delta_{42}(n)	\\
			&   &	\mbox{} + 6 \cdot \delta_{84}(n)
		+ 6 \cdot \delta_{120}(n)	\\
	a_6(n)/n	& = &	4 \cdot \delta_{30}(n)
		- 4 \cdot \delta_{60}(n)	\\
	a_7(n)/n	& = &	 \delta_{30}(n) + 4 \cdot \delta_{60}(n)
\end{eqnarray*}
hold for all the $n$ listed in Lemma~\ref{determined},
so the lemma implies that they hold for all $n \ge 3$.
These formulas imply the remarks in the introduction about the maximum
number of diagonals meeting at an interior point other than the center.
Finally
\begin{eqnarray*}
	I(n)	& = &	\delta_2(n) + \sum_{k=2}^\infty a_k(n)	\\
		& = &	\delta_2(n) + \sum_{k=2}^7 a_k(n),	\\
\end{eqnarray*}
which gives the desired formula.
(The $\delta_2(n)$ in the expression for $I(n)$ is to account
for the center point when $n$ is even, which is the only point
not counted by the $a_k$.)
\end{pf*}

%****************************************************************************
\section{The formula for regions}
\label{regions}

We now use the knowledge obtained in the proof of
Theorem~\ref{countintersections} about the number of interior points
through which exactly $k$ diagonals pass to calculate the number of
regions formed by the diagonals.

\begin{pf*}{Proof of Theorem~\ref{countregions}}
Consider the graph formed from the configuration of a regular $n$-gon
with its diagonals, in which the vertices are the vertices of the $n$-gon
together with the interior intersection points, and the edges are the
sides of the $n$-gon together with the segments that the diagonals cut
themselves into.
As usual, let $V$ denote the number of vertices of the graph, $E$ the
number of edges, and $F$ the number of regions formed, including the
region outside the $n$-gon.
We will employ Euler's Formula $V-E+F=2$.

Clearly $V=n+I(n)$.
We will count edges by counting their ends, which are $2E$ in number.
Each vertex has $n-1$ edge ends, the center (if $n$ is even) has $n$
edge ends, and any other interior point through which exactly $k$ diagonals
pass has $2k$ edge ends, so
	$$2E = n(n-1) + n \delta_2(n) + \sum_{k=2}^\infty 2k a_k(n).$$
So the desired number of regions, not counting the region outside
the $n$-gon, is
\begin{eqnarray*}
	F-1	& = & E-V+1	\\
		& = & \left[ n(n-1)/2 + n \delta_2(n)/2 +
			\sum_{k=2}^\infty k a_k(n) \right]
			- \left[ n+I(n) \right] + 1.
\end{eqnarray*}
Substitution of the formulas derived in the proof of
Theorem~\ref{countintersections} for $a_k(n)$ and $I(n)$
yields the desired result.
\end{pf*}

%****************************************************************************
\section*{Appendix: computations and tables}
\label{computations}

In Table~\ref{intersection_pts1}
we list $I(n), R(n), a_{2}(n),\ldots,a_{7}(n)$
for $n=4,5,\ldots,30$.
To determine the polynomials listed in
Theorem~\ref{countintersections}
more data was needed especially for $n \equiv 0 \bmod{6}$.
The largest
$n$ for which this was required was~420. 
For speed and memory conservation,
we took advantage of the regular $n$-gon's rotational symmetry and focused
our attention on 
only $2\pi/n$ radians of the $n$-gon. The data from this computation is
found in Table~\ref{intersection_pts2}. Although we only needed
to know the values at those $n$
listed in Lemma~\ref{determined} of Section~\ref{points},
we give a list for $n=6,12,\ldots,420$
so that the nice patterns can be seen.

The numbers in these tables were found by numerically computing (using a C
program and 64 bit precision) all possible $\left( {n}\atop{4} \right)$
intersections, and sorting them by their $x$ coordinate. We then
focused on runs of points with close $x$ coordinates, looking
for points with close $y$ coordinates. 

Several checks were made to eliminate any fears (arising from round-off
errors) of distinct points being mistaken as close. First, the C program
sent data to Maple which
checked that the coordinates of close points agreed to at least
40 decimal places. Second, we verified for each $n$ 
that close points came in counts of the
form $\left( {k}\atop{2} \right)$ ($k$ diagonals
meeting at a point give rise to $\left( {k}\atop{2} \right)$
close points. Hence, any run whose length is not of this form 
indicates a computational error). 

A second program was then written and run on a second machine to make
the computations completely rigorous. It also found the intersection
points numerically, sorted them and looked for close points,
but, to be absolutely sure that a pair of close
points $p_1$ and $p_2$ were actually the same,
it checked that for the two pairs of diagonals $(l_1,l_2)$
and $(l_3,l_4)$ determining $p_1$ and $p_2$, respectively,
the triples $l_1,l_2,l_3$ and $l_1,l_2,l_4$ each divided the
circle into arcs of lengths consistent with Theorem~\ref{thesolutions}.
Since this test only involves comparing rational numbers,
it could be performed exactly.

A word should also be said concerning limiting the search to
$2\pi/n$ radians of the $n$-gon. Both programs looked
at slightly smaller slices of the $n$-gon to avoid problems
caused by points near the boundary. We further subdivided this
region into twenty smaller pieces to make the task of sorting the
intersection points manageable. More precisely,
we limited our search to points whose angle with the
origin fell between $[c_1+ 2\pi (m-1)/(20n)+\varepsilon, 
c_1+ 2\pi m/(20n)-\varepsilon)$, $m=1,2,\ldots 20$, and also made sure
not to include the origin in the count.
Here $\varepsilon $ was chosen to be $.00000000001$
and $c_1$ was chosen to be $.00000123$ ($c_1=0$ would
have led to problems since there are many intersection points
with angle $0$ or $2\pi/n$). To make sure that
no intersection points were omitted, the number
of points found (counting multiplicity) was compared with
$({n \choose 4} -  {{n/2} \choose 2} \delta_2)/n$.

\begin{table}
\begin{center}
\begin{tabular}{|l||l|l|l|l|l|l|l|l|}
$n$&$a_{2}(n)$&$a_{3}(n)$&$a_{4}(n)$&$a_{5}(n)$&$a_{6}(n)$&$a_{7}(n)$&
$I(n)$&$R(n)$ \\ \hline \hline
3 &  &	& & & &	& 0 &				  1\\
4 &  &	& & & &	& 1 &				  4\\
5 & 5 &	& & & &	& 5 &				  11\\
6 & 12 & & & & & &13 &			          24\\
7 & 35 & & & & & & 35 &				  50\\
8 & 40 & 8 & & & & & 49	&			  80\\
9 & 126	& & & &	& & 126	&			  154\\
10 & 140 & 20 &	& & & &	161 &			  220\\
11 & 330 & & & & & & 330 &			  375\\
12 & 228 & 60 &	12 &  &	& & 301	&		  444\\
13 & 715 & & & & & & 715 &			  781\\
14 & 644 & 112 & & & & & 757 &			  952\\
15 & 1365 & & &	& & & 1365 &			  1456\\
16 & 1168 & 208	& & & &	& 1377 &		  1696\\
17 & 2380 & & &	& & & 2380 &			  2500\\
18 & 1512 & 216	& 54 & 54 & & &	1837 &		  2466\\
19 & 3876 & & &	& & & 3876 &			  4029\\
20 & 3360 & 480	& & & &	& 3841 &		  4500\\
21 & 5985 & & &	& & & 5985 &			  6175\\
22 & 5280 & 660	& & & &	& 5941 &		  6820\\
23 & 8855 & & &	& & & 8855 &			  9086\\
24 & 6144 & 864	& 264 &	24 & & & 7297 &		  9024\\
25 & 12650 & & & & & & 12650 &			  12926\\
26 & 11284 & 1196 & & &	& & 12481 &		  13988\\
27 & 17550 & & & & & & 17550 &			  17875\\
28 & 15680 & 1568 & & &	& & 17249 &		  19180\\
29 & 23751 & & & & & & 23751 &			  24129\\
30 & 13800 & 2250 & 420	& 180 &	120 & 30 & 16801& 21480\\
\hline
\end{tabular}
\end{center}
\caption{A listing of $I(n)$,$R(n)$ and $a_{2}(n),\ldots,a_{7}(n)$, 
$n=3,4,\ldots,30$.  Note that, when $n$ is even, $I(n)$ also counts
the point in the center.}
\label{intersection_pts1}
\end{table}

\begin{table}
\centerline{\small
\begin{tabular}{|l||l|l|l|l|l|l|l||l|l|l|l|l|l|l|l|}
$n$&$\frac{a_{2}(n)}{n}$&$\frac{a_{3}(n)}{n}$&$\frac{a_{4}(n)}{n}$
		&$\frac{a_{5}(n)}{n}$&
$\frac{a_{6}(n)}{n}$&$\frac{a_{7}(n)}{n}$&$\frac{I(n)-1}{n}$&
$n$&$\frac{a_{2}(n)}{n}$&$\frac{a_{3}(n)}{n}$&$\frac{a_{4}(n)}{n}$
		&$\frac{a_{5}(n)}{n}$&
$\frac{a_{6}(n)}{n}$&$\frac{a_{7}(n)}{n}$&$\frac{I(n)-1}{n}$
\\ \hline \hline
6 & 2 &   &   &   & & & 2 &                                 
	216 & 392564 & 4848 & 119 & 49 & & & 397580 \\
12 & 19 & 5 & 1 &   & & & 25 &                              
	222 & 426836 & 5166 & 126 & 54 & & & 432182 \\
18 & 84 & 12 & 3 & 3 & & & 102 &                            
	228 & 463303 & 5441 & 127 & 54 & & & 468925 \\
24 & 256 & 36 & 11 & 1 & & & 304 &                          
	234 & 501762 & 5718 & 129 & 57 & & & 507666 \\
30 & 460 & 75 & 14 & 6 & 4 & 1 & 560 &                      
	240 & 541612 & 6121 & 165 & 61 & & 5 & 547964 \\
36 & 1179 & 109 & 11 & 6 & & & 1305 &                       
	246 & 584782 & 6340 & 140 & 60 & & & 591322 \\
42 & 1786 & 194 & 27 & 13 & & & 2020 &                      
	252 & 629399 & 6693 & 137 & 70 & & & 636299 \\
48 & 3168 & 220 & 25 & 7 & & & 3420 &                       
	258 & 676580 & 6972 & 147 & 63 & & & 683762 \\
54 & 4722 & 288 & 24 & 12 & & & 5046 &                      
	264 & 725976 & 7276 & 151 & 61 & & & 733464 \\
60 & 6251 & 422 & 63 & 12 & & 5 & 6753 &                    
	270 & 777420 & 7643 & 150 & 66 & 4 & 1 & 785284 \\
66 & 9172 & 460 & 35 & 15 & & & 9682 &                      
	276 & 831575 & 7969 & 155 & 66 & & & 839765 \\
72 & 12428 & 504 & 35 & 13 & & & 12980 &                    
	282 & 887986 & 8326 & 161 & 69 & & & 896542 \\
78 & 15920 & 642 & 42 & 18 & & & 16622 &                    
	288 & 947132 & 8640 & 161 & 67 & & & 956000 \\
84 & 20007 & 805 & 43 & 28 & & & 20883 &                    
	294 & 1008358 & 9056 & 174 & 76 & & & 1017664 \\
90 & 25230 & 863 & 45 & 21 & 4 & 1 & 26164 &                
	300 & 1072171 & 9462 & 203 & 72 & & 5 & 1081913 \\
96 & 31240 & 948 & 53 & 19 & & & 32260 &                    
	306 & 1139436 & 9780 & 171 & 75 & & & 1149462 \\
102 & 37786 & 1096 & 56 & 24 & & & 38962 &                  
	312 & 1208944 & 10164 & 179 & 73 & & & 1219360 \\
108 & 45447 & 1201 & 53 & 24 & & & 46725 &                  
	318 & 1281100 & 10582 & 182 & 78 & & & 1291942 \\
114 & 53768 & 1368 & 63 & 27 & & & 55226 &                  
	324 & 1356315 & 10957 & 179 & 78 & & & 1367529 \\
120 & 62652 & 1601 & 95 & 31 & & 5 & 64384 &                
	330 & 1434110 & 11375 & 189 & 81 & 4 & 1 & 1445760 \\
126 & 73676 & 1658 & 72 & 34 & & & 75440 &                  
	336 & 1514816 & 11856 & 193 & 89 & & & 1526954 \\
132 & 85319 & 1825 & 71 & 30 & & & 87245 &                  
	342 & 1598970 & 12216 & 192 & 84 & & & 1611462 \\
138 & 97990 & 2002 & 77 & 33 & & & 100102 &                 
	348 & 1685843 & 12661 & 197 & 84 & & & 1698785 \\
144 & 112100 & 2136 & 77 & 31 & & & 114344 &                
	354 & 1775788 & 13108 & 203 & 87 & & & 1789186 \\
150 & 127070 & 2345 & 84 & 36 & 4 & 1 & 129540 &            
	360 & 1868312 & 13669 & 231 & 91 & & 5 & 1882308 \\
156 & 143635 & 2549 & 85 & 36 & & & 146305 &                
	366 & 1965272 & 14010 & 210 & 90 & & & 1979582 \\
162 & 161520 & 2736 & 87 & 39 & & & 164382 &                
	372 & 2064919 & 14465 & 211 & 90 & & & 2079685 \\
168 & 180504 & 3008 & 95 & 47 & & & 183654 &                
	378 & 2167754 & 14930 & 219 & 97 & & & 2183000 \\
174 & 201448 & 3178 & 98 & 42 & & & 204766 &                
	384 & 2274136 & 15396 & 221 & 91 & & & 2289844 \\
180 & 223251 & 3470 & 129 & 42 & & 5 & 226897 &             
	390 & 2383690 & 15885 & 224 & 96 & 4 & 1 & 2399900 \\
186 & 247562 & 3630 & 105 & 45 & & & 251342 &               
	396 & 2496999 & 16369 & 221 & 96 & & & 2513685 \\
192 & 273144 & 3844 & 109 & 43 & & & 277140 &               
	402 & 2613536 & 16896 & 231 & 99 & & & 2630762 \\
198 & 300294 & 4092 & 108 & 48 & & & 304542 &               
	408 & 2733888 & 17380 & 235 & 97 & & & 2751600 \\
204 & 329171 & 4357 & 113 & 48 & & & 333689 &               
	414 & 2857752 & 17898 & 234 & 102 & & & 2875986 \\
210 & 359556 & 4661 & 125 & 55 & 4 & 1 & 364402 &           
	420 & 2984383 & 18598 & 273 & 112 & & 5 & 3003371 \\
\hline
\end{tabular}
}
\vskip12pt
\caption{The number of intersection points for one piece
of the pie (i.e. $2\pi/n$ radians), $n=6,12,\ldots,420$.}
\label{intersection_pts2}
\end{table}

%****************************************************************************
\section*{Acknowledgements}

We thank Joel Spencer and Noga Alon for helpful conversations.
Also we thank Jerry Alexanderson, Jeff Lagarias, Hendrik Lenstra,
and Gerry Myerson for pointing out to us many of the references below.

%****************************************************************************

\end{document}